\documentclass[10pt,twoside]{article}
\usepackage{amssymb,amsmath,mathrsfs,cite}
\usepackage{latexsym}

\newcommand{\R}{{\mathbb R}}
\newcommand{\eps}{\varepsilon}
\newcommand{\de}{\partial}
\newcommand{\GG}{{\mathbb G}}
\newcommand{\Rn}{{\mathbb{R}^n}}

\newcommand{\irn}{\int_{\Rn}}
\renewcommand{\a }{\alpha }

\renewcommand{\b }{\beta }

\newcommand{\vfi}{\varphi}

\newcommand{\n }{\nabla }

\renewcommand{\t}{\theta}

\newcommand{\G}{\Gamma}

\renewcommand{\S}{\Sigma}

\renewcommand{\H}{{\mathscr H}}

\newcommand{\Ne}{{\mathcal N}}

\newtheorem{theorem}{Theorem}[section]
\newtheorem{definition}{Definition}[section]
\newtheorem{lemma}{Lemma}[section]
\newtheorem{proposition}{Proposition}[section]
\newtheorem{remark}{Remark}[section]

\newtheorem{corollary}{Corollary}[section]

\makeatletter\def\theequation{\arabic{section}.\arabic{equation}}\makeatother
\begin{document}
\title{
\vspace{0.5in}
{\bf\Large  Locating the peaks of semilinear elliptic systems}}
\author{{\bf\large Alessio Pomponio}\footnote{The author was partially supported by the MIUR research
project ``Metodi Variazionali ed Equazioni Differenziali Nonlineari''.}\hspace{2mm}\vspace{1mm}\\
{\it\small Dipartimento di Matematica} \\ {\it\small Politecnico di Bari},
{\it\small Via Amendola 126/B, I-70126 Bari, Italy}\\
{\it\small e-mail: a.pomponio@poliba.it}\vspace{1mm}\\
{\bf\large Marco Squassina}\footnote{The author was partially supported by the MIUR research
project ``Variational and Topological Methods in the Study of
Nonlinear Phenomena'' and by the Istituto Nazionale di Alta
Matematica.}\hspace{2mm}\vspace{1mm}\\
{\it\small Dipartimento di Matematica} \\ {\it\small Politecnico di Milano},
{\it\small Via Bonardi 9, I-20133 Milano, Italy}\\
{\it\small e-mail: squassina@mate.polimi.it}\vspace{1mm}}

\maketitle

\begin{center}
{\bf\small Abstract}

\vspace{3mm}
\hspace{.05in}\parbox{4.5in}
{{\small We consider a system of weakly coupled singularly perturbed
semilinear elliptic equations. First, we obtain a Lipschitz regularity
result for the associated ground energy function $\Sigma$ as
well as representation formulas for the left and the right
derivatives. Then, we show that the concentration points of
the solutions locate close to the critical points of $\Sigma$ in the sense of
subdifferential calculus.}}
\end{center}



\bigskip
\begin{center}
\begin{minipage}{12cm}
\footnotesize
\tableofcontents
\end{minipage}
\end{center}


\section{\bf Introduction and main results}
\def\theequation{1.\arabic{equation}}\makeatother
\setcounter{equation}{0}

In the asymptotic analysis of the singularly perturbed elliptic equation
\begin{equation}
\tag{$P_\eps$}
-\eps^2\Delta u+u=f(x,u) \quad\text{in $\R^n$},\qquad u>0\quad\text{in $\R^n$},
\end{equation}
there are well known situations where the associated
{\em ground energy function} $\Sigma$ (cf.\ \cite{wangzeng}) is $C^1$-smooth and
around its nondegenerate critical points the solutions $u_\eps$ of $P_\eps$
exhibit a {\em spike-like profile} as $\eps$ goes to zero.
This is the case, for instance, for the power nonlinearity
$$
f(x,u)=K(x)u^q,\qquad 1<q<{\textstyle\frac{n+2}{n-2}},\quad n\geq 3,
$$
where $K(x)$ is a suitable $C^1$ function (see e.g.\ \cite{abc,AMS} and references therein).
It turns out that the $C^1$ (and higher) smoothness of $\Sigma$
is related to the crucial fact that, for every fixed
$z\in\R^n$, the limiting autonomous equation
\begin{equation}
\tag{$P_0$}
-\Delta u+u=f(z,u) \quad\text{in $\R^n$},\qquad u>0\quad\text{in $\R^n$},
\end{equation}
admits a {\em unique solution}, up to translations \cite{wangzeng}.
However, unfortunately, the uniqueness feature for $P_0$
is a delicate matter and it is currently available only under
rather restrictive assumptions on $f$ (cf.\ e.g.~\cite{K}).
What it is know, in general, is only that $\S$ is a {\em locally Lipschitz continuous}
function which admits representation formulas for the left
and right derivatives (cf.\ \cite[Lemma 2.3]{wangzeng}).
Motivated by these facts, recently, some conditions for locating the
concentration points for $P_\eps$ in presence of a more general nonlinearity $f$, not necessarily of power type,
have been investigated (see~\cite{secsqu} and also~\cite{secsqu2}). The underlying philosophy is that
when the limit problem $P_0$ lacks of uniqueness up to translations, then the ground energy
function $\Sigma$ could loose its additional regularity properties.

Nevertheless, in this (possibly nonsmooth) framework, it turns out that a
necessary condition for the solutions $u_\eps$ to
concentrate (in a suitable sense) around a
given point $z$ is that it is critical for $\Sigma$ in the sense of the {\em Clarke
subdifferential} $\partial_C$, that is $0\in \partial_C\Sigma(z)$,
or in a even weaker sense. The main theme of this note is the search of
suitable conditions for locating the spikes, as $\eps\to 0$,
of the solutions to the semilinear model system
\begin{equation}
\label{problema}
\tag{$S_\eps$}
\begin{cases}
-\eps^2\Delta u+u=K(x)v^q, & \,\,\text{in $\R^n$}, \\
\noalign{\vskip2pt}
-\eps^2\Delta v+v=Q(x)u^p, & \,\,\text{in $\R^n$}, \\
\noalign{\vskip2pt}
\,\,\, u,\,v>0, &  \,\,\text{in $\R^n$},
\end{cases}
\end{equation}
where $p>1$ and $q>1$ are lying {\em below} the so
called ``critical hyperbola''
\[
{\mathscr C}_n=\big\{(p,q)\in (1,\infty)\times(1,\infty):\,
{\textstyle\frac{1}{p+1}}+{\textstyle\frac{1}{q+1}}=1 -{\textstyle\frac{2}{n}}\big\},
\qquad n\geq 3,
\]
which naturally arises in the study of this problem and constitutes
the borderline between existence and nonexistence results
(cf.\ e.g.\ \cite{CFM,HvdV}).

Now, according to what discussed right above,
the interest in looking for conditions for the spike location
of the solutions to \eqref{problema} is mainly motivated by the following simple observation:
contrary to the scalar case, there is {\em no uniqueness
result} available in the literature for the (radial) solutions to the
(limiting) system associated with \eqref{problema}
\begin{equation}
\label{limit-z}
\tag{$S_z$}
\begin{cases}
-\Delta u+u=K(z)v^q, & \,\,\text{in $\R^n$}, \\
\noalign{\vskip2pt}
-\Delta v+v=Q(z)u^p, & \,\,\text{in $\R^n$}, \\
\noalign{\vskip2pt}
\,\, u,\, v>0, &  \,\,\text{in $\R^n$},
\end{cases}
\end{equation}
where $z\in\R^n$ is frozen and acts as a parameter.
As a consequence, in the vectorial case,
we do not know whether the (suitably defined) ground energy map $\Sigma$
associated with \eqref{problema} (cf.\ Definition~\ref{dual-gef}) is $C^1$-{\em smooth} and
admits an {\em explicit representation formula}. Hence, the necessary
conditions in terms of Clarke subdifferential
(or weaker) appear here even more natural than in the case of a single equation.
See Section \ref{plemMa} for the statements of the main results,
Theorems \ref{neceMAIN} and \ref{neceMAIN2}.
As far as we are aware, other criteria for the concentration have been established
so far, but all of them consider the scalar case.
We refer the reader e.g.\ to~\cite{abc,pompsec,wang,wangzeng} for the
case of power-like nonlinearities and to~\cite{secsqu,secsqu2} for more
general classes of nonlinearities.

Semilinear systems like \eqref{problema} naturally arise in the
study of various kinds of nonlinear phenomena such as population evolution,
pattern formation, chemical reaction, etc., being $u$ and $v$
the concentrations of different species in the process (see also
\cite{zou} and references therein). Visibly, the interest in the study of
the various qualitative properties of \eqref{problema} has steadily increased
in recent times.
In a smooth {\em bounded} domain $\Omega$,
$(S_1)$ was extensively studied by Clement, Costa, De Figueiredo, Felmer,
Hulshof, Magalh\~aes, van der Vorst in \cite{dFM,CFM,CM,dFF,HvdV}.
The asymptotic analysis with respect to $\eps$ has been very recently performed both with
{\em Dirichlet} and {\em Neumann} boundary
conditions by Pistoia-Ramos \cite{PR,PR2} and Ramos-Yang \cite{ramyang}.
In the {\em whole space} $\R^n$, the existence of least energy solutions
to \eqref{problema} has been investigated by Alves-Carri\~ao-Miyagaki,
De Figueiredo, Yang and Sirakov in \cite{SerZou,ACM,ASY,sirak,dFJ,Jian},
whereas the asymptotic behavior with respect to $\eps$ has been pursued
by Alves-Soares-Yang in \cite{ASY}. Finally, for the exponential decay,
the radial symmetry and the regularity properties of
the solutions to \eqref{limit-z}, we refer the reader
to the quite recent achievements of Busca-Sirakov and Sirakov \cite{sirak,busi}.

The outline of the paper is as follows: in Sections \ref{dualF}-\ref{plemMa}
we provide preliminary stuff such as the (dual) variational framework and
the (dual) ground energy function $\S$ and we state the main results of the paper.
Throughout Section \ref{GeF} we deal with the ${\rm Lip}_{\rm loc}$ regularity and the representation
formulas of the directional derivatives for $\S$. Finally, in Section \ref{PROOFS}
we end up the proofs of the main results.

\subsection{The dual variational functional}
\label{dualF}
As it is known, if e.g.\ $p$ and $q$ are both less than
${\textstyle\frac{n+2}{n-2}}$, then system \eqref{problema} admits a
{\em natural variational structure} (of Hamiltonian type) which is based on
the strongly indefinite functional $f_\eps:H^1(\R^n)\times H^1(\R^n)\to\R$,
$$
f_\eps(u,v)=\int_{\R^n}\eps^2\nabla u\cdot\nabla v+uv
-{\textstyle\frac{1}{q+1}} \int_\Rn K(x)|v|^{q+1}
-{\textstyle\frac{1}{p+1}}\int_\Rn Q(x)|u|^{p+1}.
$$
However, as already done in \cite{ACM,ASY},
for our purposes, as well as for dealing with possibly supercritical
values of $p$ or $q$, we consider a corresponding
{\em dual variational structure}, mainly relying on the Legendre-Fenchel transformation
(see e.g.\ \cite{MaWi,clk,clkeke} and references therein).
In the following, we just briefly recall some of the core ingredients,
referring to \cite[Section 2]{ACM} for expanded details on this framework.
For ${\textstyle\frac{1}{p+1}}+{\textstyle\frac{1}{q+1}}>{\textstyle\frac{n-2}{n}}$,
consider the linear operators
\begin{align*}
& T_1:L^{\frac{q+1}{q}}(\R^n)\to W^{2,\frac{q+1}{q}}(\R^n)\hookrightarrow L^{p+1}(\R^n), \\
& T_2:L^{\frac{p+1}{p}}(\R^n)\to W^{2,\frac{p+1}{p}}(\R^n)\hookrightarrow L^{q+1}(\R^n),
\end{align*}
defined as
$$
T_1=T_2=(-\Delta+{\rm Id})^{-1}.
$$
Notice that $T_1$ and $T_2$ are continuous.
Then, we consider the linear operator (take into account the proper
Sobolev embeddings)
$$
T:L^{\frac{p+1}{p}}(\R^n)\times L^{\frac{q+1}{q}}(\R^n)\to
L^{p+1}(\R^n)\times L^{q+1}(\R^n),
\qquad
T=\begin{bmatrix}
0 & T_1 \\
T_2& 0
\end{bmatrix},
$$
explicitly defined by
$$
\langle T\eta,\xi\rangle=
\xi_1T_1\eta_2+\xi_2T_2\eta_1,
\qquad\forall \eta=(\eta_1,\eta_2),\,\,\,\forall\xi=(\xi_1,\xi_2).
$$
Finally we introduce the Banach space $({\H},\|\cdot\|_{\H})$,
$$
{\H}=L^{\frac{p+1}{p}}(\R^n)\times L^{\frac{q+1}{q}}(\R^n),
\qquad
\|\eta\|_{\H}^2=\|\eta_1\|^2_{L^{\frac{p+1}{p}}(\R^n)}
+\|\eta_2\|^2_{L^{\frac{q+1}{q}}(\R^n)}
$$
and the (dual) $C^1$ functional $J_\eps:{\H}\to\R$ defined as
$$
J_\eps(\eta)={\textstyle\frac{p}{p+1}}\int_{\R^n}\frac{|\eta_1|^\frac{p+1}{p}}{Q^{\frac{1}{p}}(\eps x)}+
{\textstyle\frac{q}{q+1}}\int_{\R^n}\frac{|\eta_2|^\frac{q+1}{q}}{K^{\frac{1}{q}}(\eps x)}-
{\textstyle\frac{1}{2}}\int_{\R^n}\langle T\eta,\eta\rangle .
$$
If $\eta^\eps=(\eta^\eps_1,\eta^\eps_2)$ is a critical point of $J_\eps$,
then $(u_\eps(x),v_\eps(x))=(\bar u_\eps(\frac{x}{\eps}),\bar v_\eps(\frac{x}{\eps}))$, with
\begin{equation}
\label{summab}
(\bar u_\eps,\bar v_\eps)=(T_1\eta^\eps_2,T_2\eta^\eps_1)\in W^{2,\frac{q+1}{q}}\cap L^{p+1}\times
W^{2,\frac{p+1}{p}} \cap L^{q+1},
\end{equation}
corresponds to a solution to \eqref{problema} with $u_\eps(x),v_\eps(x)\to 0$
for $|x|\to\infty$ (see \cite[p.677]{ACM}). In light of
the above summability, we have $f_\eps(u_\eps,v_\eps)\in\R$ for all $\eps>0$.
Analogously, associated with \eqref{limit-z}, we introduce the limiting
functional $I_z:{\H}\to\R$
\begin{equation*}
I_z(\eta)={\textstyle\frac{p}{p+1}}\int_{\R^n}\frac{|\eta_1|^\frac{p+1}{p}}{Q^{\frac{1}{p}}(z)}+
{\textstyle\frac{q}{q+1}}\int_{\R^n}\frac{|\eta_2|^\frac{q+1}{q}}{K^{\frac{1}{q}}(z)}-
{\textstyle\frac{1}{2}}\int_{\R^n}\langle T\eta,\eta\rangle .
\end{equation*}
From the viewpoint of our investigation, the main advantage of exploiting
the dual variational functional $I_z$ is that
it admits a {\em mountain-pass geometry} and the mountain-pass
value corresponds to the {\em least possible energy}
of system \eqref{limit-z}. As we shall see in the next section,
this allows to provide in the vectorial framework a suitable
definition of ground energy function with nice features, similar to those available in
the scalar case.

\subsection{Preliminaries and the main results}
\label{plemMa}

In order to state the main achievements of the paper, we need
some preparatory stuff. For the sake of self-containedness
we shall also recall a few pretty well known notions from
nonsmooth calculus (see e.g.~\cite{clarke}).

\begin{definition}\label{def:der}\rm
Let $f:\R^n\to\R$ be a locally Lipschitz function near a point $z\in\R^n$.
The {\em Clarke subdifferential} of $f$ at $z$ is defined by
\begin{equation*}
\partial_C f(z):=\big\{\eta\in\R^n:\,\,
f^0(z,w)\geq \eta\cdot w,\,\,\,\text{for every $w\in\R^n$}\big\},
\end{equation*}
where $f^0(z,w)$ is the generalized derivative of
$f$ at $z$ along $w\in\R^n$, defined by
\begin{equation*}
f^0(z;w):=\limsup_{\substack{\xi \to z
\\ \lambda \to 0+}}\frac{f(\xi+\lambda w)-f(\xi)}{\lambda}.
\end{equation*}
\end{definition}

\begin{definition}\rm
\label{dual-gef}
The (dual) {\em ground energy function} $\Sigma:\R^n\to\R$ of~\eqref{limit-z} is given by
\begin{equation*}
\Sigma(z):=\inf_{\eta \in {\mathcal N}_z}I_z (\eta),
\end{equation*}
where ${\mathcal N}_z$ is the {\em Nehari manifold} of $I_z$, that is
\[
{\mathcal N}_z = \big\{\eta \in {\H}:\,\,
\text{$\eta\neq (0,0)$ and $I_z' (\eta)[\eta]=0$}
\big\}.
\]
We shall denote by $\mathcal{K}\subset\R^n$ the set of
{\em Clarke critical points} of $\Sigma$, namely
$$
\mathcal{K}:=\big\{z\in\R^n:\, 0\in\partial_C\Sigma(z)\big\}.
$$
\end{definition}

\begin{definition}\rm
We say that the pair $(u_{\eps},v_{\eps})$ is a {\em strong solution} to system~\eqref{problema}
if it is a distributional solution and $(u_{\eps},v_{\eps})\in
W^{2,(q+1)/q}(\R^n)\times W^{2,(p+1)/p}(\R^n)$. We say that
the pair $\eta^\eps=(\eta_1^\eps,\eta_2^\eps)$
corresponding to $(u_{\eps},v_{\eps})$ through \eqref{summab} is the
related {\em dual solution}.
\end{definition}

\begin{definition}\rm
We set
\begin{align*}
\mathcal{E}:=\big\{& z\in\R^n:\, \text{there exists a sequence of strong solutions
$(u_{\eps_h},v_{\eps_h})$ of~\eqref{problema} with} \\
& \text{$|u_{\eps_h}(z)|,|v_{\eps_h}(z)|\geq\delta$
for some $\delta>0$, $|u_{\eps_h}(z+\eps_hx)|,\,|v_{\eps_h}(z+\eps_hx)|\to 0$} \\
&\text{as $|x|\to\infty$ uniformly w.r.t.\ $h$, and
${\eps_h}^{-n}f_{\eps_h}(u_{\eps_h},v_{\eps_h})\to\Sigma(z)$ as $h\to\infty$}
\big\}.
\end{align*}
We say that $\mathcal{E}$ is the {\em energy concentration set} for~\eqref{problema}.
\end{definition}

Assume that $K,Q\in C^1(\R^n)$ and
\begin{equation}
\label{asspoteorig}
\a \leq  K(x)\leq \b,
\quad
\a \leq  Q(x)\leq \b,
\qquad\text{for all $x\in\R^n$},
\end{equation}
\begin{equation}
\label{assp2}
|\nabla K(x)|,
|\nabla Q(x)|\leq Ce^{M|x|},
\qquad\text{for all $x\in\R^n$ with $|x|$ large}.
\end{equation}
for some positive constants $\alpha,\b,C$ and $M$.
\vskip4pt

The main result of the paper, linking the energy
concentration set $\mathcal{E}$ with the set $\mathcal{K}$
of Clarke critical set of $\Sigma$, is provided by the following

\begin{theorem}
\label{neceMAIN}
Assume that $K,Q\in C^1(\R^n)$ and that
\eqref{asspoteorig}-\eqref{assp2} hold.
Then $\mathcal{E}\subset\mathcal{K}$.
\end{theorem}

\begin{remark}\rm
By \cite[Theorem 1]{ASY}, under suitable assumptions, if there exists an absolute
minimum (or maximum) point $z_*$ for $\Sigma$, then $z_*\in\mathcal{E}\not=\emptyset$.
\end{remark}

\begin{remark}\rm
As a straightforward combination of Theorem \ref{neceMAIN} with the
well known convex hull characterization of $\partial_C\Sigma(z)$, if $z$ is a concentration
point for~\eqref{problema}, then
$$
0\in {\rm Co}\big\{\lim_j\nabla\Sigma(\xi_j):
\,\text{$\xi_j\not\in\Omega$ and $\xi_j\to z$}\big\},
$$
where ${\rm Co}\{X\}$ denotes the convex hull of $X$ and $\Omega$ is any
null set containing the set of points at which $\Sigma$ fails to be
differentiable.
\end{remark}

\begin{corollary}
\label{uniTH}
Under the (unproved) assumption that, for all $z\in\R^n$, system \eqref{limit-z} admits a unique
positive solution (up to translations), $\S$ is $C^1$-smooth and
$$
\mathcal{E}\subset{\rm Crit}\big(Q^{\frac{q+1}{pq-1}}\, K^{\frac{p+1}{pq-1}}\big),
$$
where {\rm Crit}(f) denotes the set of (classical) critical points of $f$.
\end{corollary}

In the following definition we consider solutions which concentrate
close to a point $z$, with bounded energy but not necessary stabilizing towards $\S(z)$.

\begin{definition}\rm
Let $m\geq 1$. We set
\begin{align*}
\mathcal{E}_m:=\big\{& z\in\R^n:\, \text{there exists a sequence of strong solutions
$(u_{\eps_h},v_{\eps_h})$ of~\eqref{problema} with} \\
& \text{$|u_{\eps_h}(z)|,|v_{\eps_h}(z)|\geq\delta$
for some $\delta>0$, $|u_{\eps_h}(z+\eps_hx)|,\,|v_{\eps_h}(z+\eps_hx)|\to 0$} \\
&\text{as $|x|\to\infty$ uniformly w.r.t.\ $h$, and
${\eps_h}^{-n}f_{\eps_h}(u_{\eps_h},v_{\eps_h})\to m$ as $h\to\infty$}
\big\}.
\end{align*}
We say that $\mathcal{E}_m$ is the {\em concentration set} for~\eqref{problema}
at the energy level $m$.
\end{definition}

\begin{definition}\rm
\label{gammapm}
Let $m\geq 1$ and $z\in\R^n$. For every $w\in\R^n$ we define $\Gamma_{z,m}^\mp(w)$ by
\begin{align*}
\Gamma_{z,m}^-(w)&:=\sup_{\eta \in {\GG}_m(z)}\bigg[
-\frac{1}{p+1}\frac{\de Q}{\de w}(z) \int_{\R^n}\frac{|\eta_1|^\frac{p+1}{p}}{Q^{\frac{p+1}{p}}(z)}
-\frac{1}{q+1}\frac{\de K}{\de w}(z) \int_{\R^n}\frac{|\eta_2|^\frac{q+1}{q}}{K^{\frac{q+1}{q}}(z)}\bigg], \\
\noalign{\vskip4pt}
\Gamma_{z,m}^+(w)&:=-\inf_{\eta \in {\GG}_m(z)}\bigg[
-\frac{1}{p+1}\frac{\de Q}{\de w}(z) \int_{\R^n}\frac{|\eta_1|^\frac{p+1}{p}}{Q^{\frac{p+1}{p}}(z)}
-\frac{1}{q+1}\frac{\de K}{\de w}(z) \int_{\R^n}\frac{|\eta_2|^\frac{q+1}{q}}{K^{\frac{q+1}{q}}(z)}
\bigg],
\end{align*}
where $\GG_m(z)$ denotes the set of all the nontrivial, radial, exponentially decaying
solutions of \eqref{limit-z} having energy equal to $m$.
\end{definition}

It is readily seen that $\Gamma_{z,m}^\mp(w)\in\R$
for all $z,w$ in $\R^n$ (see the proof of \eqref{Lbdd-duale}).
It is also straightforward to check that, for any $z\in\R^n$, the functions
$\{w\mapsto \Gamma_{z,m}^\mp(w)\}$ are {\em convex}.

\begin{definition}\rm
Let $m\geq 1$. We set
$$
\mathcal{K}_m:=\Big\{z\in\R^n:\,\,0\in\partial\Gamma_{z,m}^-(0)\cap\partial\Gamma_{z,m}^+(0)\Big\},
$$
where $\partial$ stands for the subdifferential of convex functions,
$$
\partial\Gamma_{z,m}^\mp(0)=\Big\{\xi\in\R^n:\,\,
\text{$\Gamma_{z,m}^\mp(w)\geq \xi\cdot w$,\,\,
for every $w\in\R^n$}\Big\}.
$$
\end{definition}

It is known by standard convex analysis that $\partial\Gamma_{z,m}^\mp(0)\not=\emptyset$,
for every $z\in\R^n$. Observe that $z\in\mathcal{K}_m$ if and only if $0$ is a critical
point for both $\Gamma_{z,m}^-$ and $\Gamma_{z,m}^+$. Of course, if $\GG_m(z)=\{\eta_0\}$
was a singleton, then $z\in\mathcal{K}_m$ if and only if
$$
\Gamma_{z,m}^-(w)=\Gamma_{z,m}^+(w)=\frac{\partial I_z}{\partial w}(\eta_0)=0,\qquad\forall w\in\R^n.
$$
\vskip2pt
Without forcing the energy levels of the solutions to approach
the least energy of the limit system, we get the following correlation
between the sets $\mathcal{E}_m$ and $\mathcal{K}_m$.

\begin{theorem}
\label{neceMAIN2}
Assume that $K,Q\in C^1(\R^n)$ and
\eqref{asspoteorig}-\eqref{assp2} hold.
Then $\mathcal{E}_m\subset\mathcal{K}_m$.
\end{theorem}

\section{\bf Properties of the ground energy function}
\label{GeF}
\def\theequation{2.\arabic{equation}}\makeatother
\setcounter{equation}{0}

Before coming to the proof of the results, we need
some preliminary stuff.

\subsection{Some preparatory lemmas}
The next proposition is well known (see e.g. \cite{Jian}); on the other hand,
for the sake of completeness and  self-containedness, we report a brief proof.

\begin{proposition}
\label{equiv-sol}
Let $z\in\R^n$. Then $(u,v)\in W^{2,\frac{q+1}{q}}(\R^n)\times W^{2,\frac{p+1}{p}}(\R^n)$ is a
solution to \eqref{limit-z} if and only if $\eta=(\eta_1,\eta_2)=(T^{-1}_2 v,T^{-1}_1 u)
\in\H$ is a critical point of $I_z$. Moreover, there holds
$f_z(u,v)=I_z(\eta_1,\eta_2)$, where $f_z$ is the functional defined as
$$
f_z(u,v)=\int_{\R^n}\nabla u\cdot\nabla v+uv
-{\textstyle\frac{1}{q+1}} \int_\Rn K(z)|v|^{q+1}
-{\textstyle\frac{1}{p+1}}\int_\Rn Q(z)|u|^{p+1}.
$$
\end{proposition}
\begin{proof}
Observe first that, if $(u,v)\in W^{2,\frac{q+1}{q}}(\R^n)
\times W^{2,\frac{p+1}{p}}(\R^n)$ solves \eqref{limit-z},
taking into account the Sobolev embedding, the value $f_z(u,v)$
is indeed finite (cf.\ \eqref{summab}). Let $(u,v)$ be a
solution to \eqref{limit-z}. Then, since
\[
\eta_1=T^{-1}_2 v, \qquad \eta_2=T^{-1}_1 u,
\]
we have
\begin{equation*}
\begin{cases}
\eta_2=T^{-1}_1 u= -\Delta u+u=K(z)v^q, \\
\noalign{\vskip2pt}
\eta_1=T^{-1}_2 v=-\Delta v+v=Q(z)u^p.
\end{cases}
\end{equation*}
Therefore, we get
\begin{equation}
\label{Tjid}
T_2 \eta_1 =v= \frac{\eta_2^{\frac 1q}}{K^{\frac 1q}(z)}
\qquad {\rm and} \qquad
T_1 \eta_2 =u= \frac{\eta_1^{\frac 1p}}{Q^{\frac 1p}(z)},
\end{equation}
and so $(\eta_1, \eta_2)$ is a critical point of $I_z$.
Vice versa, if $(\eta_1, \eta_2)$ is a critical point of $I_z$,
it is readily seen that \eqref{Tjid} hold, so that $(T_1\eta_2,T_2\eta_1)=(u,v)$
is a solution to \eqref{limit-z} (cf.\ \cite[p.677]{ACM}).
Furthermore, on the solutions to \eqref{limit-z}, we have
\[
f_z(u,v)=
\big({\textstyle\frac 12 - \frac{1}{p+1}}\big) \int_\Rn Q(z) u^{p+1}+
\big({\textstyle\frac 12 - \frac{1}{q+1}}\big) \int_\Rn K(z) v^{q+1}.
\]
Then, in light of \eqref{Tjid}, we have
\begin{align*}
I_z(\eta)&={\textstyle\frac{p}{p+1}}\int_{\R^n}\frac{|\eta_1|^\frac{p+1}{p}}{Q^{\frac{1}{p}}(z)}+
{\textstyle\frac{q}{q+1}}\int_{\R^n}\frac{|\eta_2|^\frac{q+1}{q}}{K^{\frac{1}{q}}(z)}-
{\textstyle\frac{1}{2}}\int_{\R^n}\langle T\eta,\eta\rangle  \\
\noalign{\vskip2pt}
&= \big( {\textstyle\frac{p}{p+1} - \frac 12 }\big) \irn \eta_1 \, T_1 \eta_2
+\big( {\textstyle\frac{q}{q+1} - \frac 12 }\big) \irn \eta_2 \, T_2 \eta_1 \\
\noalign{\vskip2pt}
&= \big({\textstyle \frac{p}{p+1} - \frac 12 }\big) \irn (- \Delta v +v) u
+\big( {\textstyle\frac{q}{q+1} - \frac 12 }\big) \irn (- \Delta u +u) v \\
\noalign{\vskip2pt}
&=\big({\textstyle\frac 12 - \frac{1}{p+1}}\big) \int_\Rn Q(z) u^{p+1}+
\big({\textstyle\frac 12 - \frac{1}{q+1}}\big) \int_\Rn K(z) v^{q+1}
=f_z(u,v),
\end{align*}
which concludes the proof.
\end{proof}

\begin{definition}\rm
We say that $\eta\in\H$ is a  {\em dual solution} to \eqref{limit-z}
if it is a critical point of $I_z$. We say that $\eta$ is a  {\em dual least energy
solution} to \eqref{limit-z} if it is a dual solution and, in addition, $I_z(\eta)=\Sigma(z)$.
\end{definition}

The next property, classical in the scalar case, will be pretty useful
for our purposes.

\begin{lemma}
For every $z\in\R^n$, let us set
\begin{align*}
b_1(z)&:=\inf_{\eta\in {\H}\setminus \{0\}} \sup_{t \ge 0} I_z(t\eta), \\
\noalign{\vskip5pt}
b_2(z)&:=\inf_{\eta \in {\mathcal N}_z}I_z (\eta)=\Sigma(z), \\
\noalign{\vskip5pt}
b_3(z)&:=\inf\big\{I_z (\eta):\text{$\eta
\in{\H}\setminus \{0\}$ is a dual solution to \eqref{limit-z}}\big\}.
\end{align*}
Then $b_1(z)=b_2(z)=b_3(z)$. Moreover $\{z\mapsto\S(z)\}$ is continuous.
\end{lemma}
\begin{proof}
The first equality follows from \cite[Lemma 2]{ACM}. Moreover in \cite{ACM} it is proved
that $b_1(z)=b_2(z)$ is a critical value so that also $b_2(z)=b_3(z)$ follows.
Finally, by virtue of \cite[Lemma 1]{ASY}, we know that $\S$ is continuous.
\end{proof}

\begin{lemma}
\label{le:retta}
Let $z\in\R^n$ and define the (nonempty) set
\[
\H_+ :=\Big\{\eta\in \H :\,\irn  \langle T \eta,\eta \rangle >0\Big\}.
\]
Then, for every $\eta \in \H_+$, there exists a unique maximum point
$t_{\eta}>0$ of the map $\phi \colon t\in (0,\infty) \mapsto I_z(t \eta)$.
In particular, $t_\eta \eta \in {\mathcal N}_z$.
\end{lemma}
\begin{proof}
Let us observe that if $\phi'(t)=0$, then
\[
\irn \langle T \eta,\eta \rangle
=t^{\frac{1-p}{p}} \irn \frac{|\eta_1|^{\frac{p+1}{p}}}{Q^{\frac 1p}(z)}
+t^{\frac{1-q}{q}} \irn \frac{|\eta_2|^{\frac{q+1}{q}}}{K^{\frac 1q}(z)} .
\]
Since the function $g(t)=A t^{\frac{1-q}{q}}+B t^{\frac{1-p}{p}}$ with $A,B>0$
is strictly decreasing for $t>0$, then
$\phi$ has at most one critical value. It is easy to see that for all $\eta\in \H$, $\phi(t)>0$ for
$t$ small, while if $\eta\in \H_+$, it is readily seen that $\phi(t)<0$ for big $t$'s.
\end{proof}

\subsection{Conjecturing the representation of $\S$}
\label{uniqrem}
Consider for a moment the equation
\begin{equation}
\label{eq:eq}
-\eps^2 \Delta u+V(x)u=K(x)u^p, \qquad \text{in $\R^n$},
\end{equation}
with $p$ subcritical and $V$ and $K$ potentials functions bounded away from zero.
By the results of~\cite{K}, we know that there is uniqueness (up to translation)
of positive solutions for
\begin{equation*}
-\Delta u+u=u^p, \qquad \text{in $\R^n$},
\end{equation*}
and, by a suitable change of variable, also for
the ``limit" problem at $x=z$ of \eqref{eq:eq}
\begin{equation*}
-\Delta u+V(z)u=K(z)u^p, \qquad \text{in $\R^n$}.
\end{equation*}
This allows to give an explicit representation for the ground state function
associated with~\eqref{eq:eq}, merely depending on the potentials $V$ and $K$
(see for example \cite{AMS,wangzeng}):
\begin{equation}
\label{sigmarap-sc}
\S(z)=\G \,\frac{V^{\frac{p+1}{p-1}- \frac n2}(z)}{K^{\frac{2}{p-1}}(z)},
\end{equation}
for a suitable positive constant $\G$.
On the contrary, as already observed, up to our knowledge
there is no (known) uniqueness result for the elliptic system
\begin{equation}
\label{eq:sistema}
-\Delta \xi + \xi = \zeta^q, \qquad
-\Delta \zeta + \zeta = \xi^p, \qquad\text{in $\R^n$},
\end{equation}
and so, in general, we cannot provide an explicit expression for $\S$ for \eqref{limit-z}.
Slightly more in general, if $V$ is smooth and $\alpha\leq V(x)\leq\beta$, consider the system
\begin{equation}
\label{limit-zbis}
\begin{cases}
-\Delta u+V(z)u=K(z)v^q, & \,\,\text{in $\R^n$}, \\
\noalign{\vskip2pt}
-\Delta v+V(z)v=Q(z)u^p, & \,\,\text{in $\R^n$}, \\
\noalign{\vskip2pt}
\,\, u,\,v>0, &  \,\,\text{in $\R^n$}.
\end{cases}
\end{equation}
{\em Assuming} for a moment that~\eqref{eq:sistema} has a
{\em unique} solution $(\xi,\zeta)$, then we claim that
\begin{equation}
\label{conjsigma}
\Sigma(z)=
\Gamma\frac{V^{\frac{(p+1)(q+1)}{pq-1}-\frac{n}{2}}(z)}
{Q^{\frac{q+1}{pq-1}}(z)\, K^{\frac{p+1}{pq-1}}(z)},
\end{equation}
for a suitable positive constant $\G$. Indeed, by rescaling
\[
u(x)=\varpi_1\xi(\mu x)
\quad
\text{and}
\quad
v(x)=\varpi_2\zeta(\mu x),
\]
where we have set
\begin{align*}
\mu=\mu(z) &:=V^\frac 12 (z), \\
\varpi_1=\varpi_1(z) &:=
\frac{V^{\frac{q+1}{pq-1}}(z)}{Q^{\frac{q}{pq-1}}(z)\,K^{\frac{1}{pq-1}}(z)}, \\
\varpi_2=\varpi_2(z) &:=
\frac{V^{\frac{p+1}{pq-1}}(z)}{Q^{\frac{1}{pq-1}}(z)\,K^{\frac{p}{pq-1}}(z)},
\end{align*}
it is easy to see that $(u, v)$ is the unique solution of the system \eqref{limit-zbis}.
Hence, by a straightforward calculation, we reach \eqref{conjsigma}.
Let us observe that the exponent of $V(z)$ in~\eqref{conjsigma} is equal to zero if, and only if,
the pair $(p,q)$ belongs to ${\mathscr C}_n$.
Then, for problems with powers $p,q$ close to the set ${\mathscr C}_n$, the potential $V$
is expected to have a weak influence in the location of concentration points.
Notice that the same phenomenon appears in the scalar case (cf.\ formula~\eqref{sigmarap-sc}),
since $\frac{p+1}{p-1}- \frac n2\sim 0$ if and only if $p\sim \frac{n+2}{n-2}=2^*-1$, where
$2^*$ is the critical Sobolev exponent for $H^1$.
Finally, we just wish to mention that, incidentally, the exponents
$$
\theta_1=\frac{p+1}{pq-1},\qquad \theta_2=\frac{q+1}{pq-1}
$$
in formula \eqref{conjsigma} also arise
in the study of the {\em blow-up rates} for the
parabolic system
$$
u_t=\Delta u+v^q, \qquad v_t=\Delta v+u^p,\qquad  x\in\Omega,\,\, t>0,
$$
with initial data $u(x,0)=u_0(x)\geq 0$, $v(x,0)=v_0(x)\geq 0$ and
Dirichlet boundary conditions $u=v=0$ on $\partial\Omega$. Here $\Omega$ is a
ball in $\R^n$ and $u_0$ and $v_0$ are continuous which vanish
on the boundary. If $u_0$, $v_0$ are nontrivial $C^1$ functions,
the solution $(u,v)$ blows up at {\em a finite time $T<\infty$}, and $u_t\ge 0$, $v_t\ge 0$
on $\Omega\times(0,T)$, then there exist two constants $C>c>0$ with
$$
\frac{c}{(T-t)^{\theta_1}}\leq \max_{\overline\Omega}u(x,t)\leq \frac{C}{(T-t)^{\theta_1}},\qquad
\frac{c}{(T-t)^{\theta_2}}\leq \max_{\overline\Omega}v(x,t)\leq \frac{C}{(T-t)^{\theta_2}},
$$
for all $t\in(0,T)$. We refer the interested reader, e.g., to \cite{mwang}.

\subsection{Local lipschitzianity of $\Sigma$}

In the case of a single semilinear elliptic equation, it is known~\cite{wangzeng}
that the ground energy map enjoys a basic regularity property, in addition
to the continuity, namely it is locally Lipschitz continuous (hence differentiable {\em a.e.}\ by
virtue of Rademacher's theorem). Analogously, for system~\eqref{problema}, we
obtain the following

\begin{theorem}
\label{loclips}
$\Sigma\in {\rm Lip}_{\rm loc}(\R^n)$.
\end{theorem}
\begin{proof}
Let $\rho_0>0$ and $\mu\in\R^n$ with $|\mu|\leq \rho_0$ and
let $\eta^\mu$ be a (dual) solution to ($S_\mu$) such that
$I_\mu(\eta^\mu)=\Sigma(\mu)$ (we already know that such a solution
does exist, see \cite{ACM}). Then, the corresponding (direct)
solution $(u_\mu,v_\mu)$ satisfies
\begin{equation}
\label{eq:Lmu}
-\Delta u + u = K(\mu) v^q, \qquad
-\Delta v + v = Q(\mu) u^p, \qquad\text{in $\R^n$}.
\end{equation}
We also know that $u_\mu$ and $v_\mu$ are radially symmetric,
radially decreasing with respect to, say, the origin, and exponentially
decaying (see \cite{dFJ,sirak,busi}, in particular \cite[Theorem 2]{busi}
and \cite[Theorem 1(a)]{sirak}). We claim that there exist
$\varpi_1>0$ and $\varpi_2=\varpi_2(\rho_0)>0$ independent of $\mu$ such that
\begin{equation}
\label{Lbdd}
\varpi_1\leq \|u_\mu\|_{L^{p+1}}\leq \varpi_2
\qquad
\text{and}
\qquad
\varpi_1\leq\|v_\mu\|_{L^{q+1}}\leq \varpi_2.
\end{equation}
Let us prove first the estimates from below.
By multiplying the first equation of~\eqref{eq:Lmu} by $u_\mu$ and
taking into account \eqref{asspoteorig}, we get
\begin{align}
\label{pezzo1}
\|u_\mu\|_{H^1}^2=\int_{\R^n}K(\mu)v^q_\mu u_\mu
\leq\b\|v_\mu\|_{L^{q+1}}^{q}\|u_\mu\|_{L^{q+1}}
\leq \b S\|v_\mu\|_{L^{q+1}}^{q}\|u_\mu\|_{H^1},
\end{align}
where $S$ is the Sobolev constant. Now, by multiplying the first equation of system~\eqref{eq:Lmu} by $v_\mu$
and the second equation by $u_\mu$, and comparing the resulting equations, we have
\begin{equation}
\label{pezzo2}
\|v_\mu\|^q_{L^{q+1}}\leq {\textstyle\left(\frac{\b}{\a}\right)^{q/(q+1)}}
\|u_\mu\|_{L^{p+1}}^{q(p+1)/(q+1)}.
\end{equation}
By combining inequalities~\eqref{pezzo1} and~\eqref{pezzo2}, and using
again the Sobolev inequality, the assertion follows.
The proof of the estimate from below for $\|v_\mu\|_{L^{q+1}}$ is similar.
To prove the inequalities from above we simply observe that $\Sigma$ is continuous and
\begin{align*}
\max_{|\mu|\leq\rho_0}\S(\mu)&=\max_{|\mu|\leq\rho_0}I_\mu(\eta^\mu)=\max_{|\mu|\leq\rho_0}f_\mu(u_\mu, v_\mu) \\
&\geq {\textstyle\big(\frac \a 2 - \frac{\a}{q+1}\big)} \|v_\mu\|_{L^{q+1}}^{q+1}
+{\textstyle\big(\frac \a 2 - \frac{\a}{p+1}\big)} \|u_\mu\|_{L^{p+1}}^{p+1}.
\end{align*}
Thus \eqref{Lbdd} follows. As a consequence, according to the definition of the dual norm
$\|\cdot\|_\H$, we immediately obtain
\begin{equation}
\label{Lbdd-duale}
\alpha\sqrt{\varpi_1^{2p}+\varpi_1^{2q}}\leq
\max_{|\mu|\leq\rho_0}\|\eta^\mu\|_{\H}\leq \beta\sqrt{\varpi_2^{2p}(\rho_0)+\varpi_2^{2q}(\rho_0)}.
\end{equation}
Now, since $\eta^\mu \in {\mathcal N}_\mu$, we get
\begin{equation}
\label{N-mu}
\irn \langle T \eta^\mu, \eta^\mu \rangle
=\int_{\R^n}\frac{|\eta_1^\mu|^\frac{p+1}{p}}{Q^{\frac{1}{p}}(\mu)}
+\int_{\R^n}\frac{|\eta_2^\mu|^\frac{q+1}{q}}{K^{\frac{1}{q}}(\mu)}>0.
\end{equation}
Hence $\eta^\mu \in \H_+$ and, by means of Lemma \ref{le:retta}, there exists precisely
one positive number $\theta (\mu,\xi)$ such that
$\theta (\mu,\xi) \eta^\mu \in{\mathcal N}_\xi$.
By definition, this means that
\begin{equation}\label{N-xi}
\irn \langle T \eta^\mu, \eta^\mu \rangle
=\theta (\mu,\xi)^{\frac{1-p}{p}}
\int_{\R^n}\frac{|\eta_1^\mu|^\frac{p+1}{p}}{Q^{\frac{1}{p}}(\xi)}
+\theta (\mu,\xi)^{\frac{1-q}{q}}
\int_{\R^n}\frac{|\eta_2^\mu|^\frac{q+1}{q}}{K^{\frac{1}{q}}(\xi)}.
\end{equation}
Moreover, we have $\theta (\mu,\mu)=1$. Collecting these facts, we see that, by the
implicit function theorem, $\t$ is differentiable with respect to the
variable $\xi$. Moreover, in light of \eqref{Lbdd-duale}, it results that
$\theta (\mu,\xi)$ remains bounded for $\mu$ and $\xi$ varying in a bounded set.
Indeed, by combining \eqref{N-mu} and \eqref{N-xi},
supposing for example that $p\le q$, we have
\[
\theta (\mu,\xi)^{\frac{p-1}{p}} \left[
 \int_{\R^n}\frac{|\eta_1^\mu|^\frac{p+1}{p}}{Q^{\frac{1}{p}}(\mu)}
+\int_{\R^n}\frac{|\eta_2^\mu|^\frac{q+1}{q}}{K^{\frac{1}{q}}(\mu)}  \right]
= \int_{\R^n}\frac{|\eta_1^\mu|^\frac{p+1}{p}}{Q^{\frac{1}{p}}(\xi)}
+\theta (\mu,\xi)^{\frac{1}{q}- \frac 1p}
\int_{\R^n}\frac{|\eta_2^\mu|^\frac{q+1}{q}}{K^{\frac{1}{q}}(\xi)}.
\]
Then the (local) boundedness of $\theta (\mu,\xi)$ follows immediately
by \eqref{Lbdd-duale} and by the fact that $\frac{1}{q}-\frac{1}{p} \le 0$.
Let us now observe that
\begin{align*}
I_\xi (\theta (\mu,\xi) \eta^\mu)
&=\theta (\mu,\xi)^{\frac{p+1}{p}} \frac{p}{p+1}
\int_{\R^n}\frac{|\eta_1^\mu|^\frac{p+1}{p}}{Q^{\frac{1}{p}}(\xi)} \\
\noalign{\vskip2pt}
&+\theta (\mu,\xi)^{\frac{q+1}{q}} \frac{q}{q+1}
\int_{\R^n}\frac{|\eta_2^\mu|^\frac{q+1}{q}}{K^{\frac{1}{q}}(\xi)} \\
\noalign{\vskip3pt}
&- \frac{\theta (\mu,\xi)^2}{2}\irn \langle T \eta^\mu, \eta^\mu \rangle .
\end{align*}
The gradient of the function
$\big\{\xi \mapsto I_\xi (\theta (\mu,\xi) \eta^\mu)\big\}$
is thus given by
\begin{align*}
\n_\xi I_\xi (\theta (\mu,\xi) \eta^\mu)\!
=& - \frac{\theta (\mu,\xi)^{\frac{p+1}{p}}}{p+1} \n_\xi Q(\xi)
\int_{\R^n}\frac{|\eta_1^\mu|^\frac{p+1}{p}}{Q^{\frac{p+1}{p}}(\xi)} \\
&- \frac{\theta (\mu,\xi)^{\frac{q+1}{q}}}{q+1} \n_\xi K(\xi)
\int_{\R^n}\frac{|\eta_2^\mu|^\frac{q+1}{q}}{K^{\frac{q+1}{q}}(\xi)}
\\
&+\n_\xi \theta (\mu,\xi)
\bigg[\theta (\mu,\xi)^{\frac{1}{p}}
\int_{\R^n}\frac{|\eta_1^\mu|^\frac{p+1}{p}}{Q^{\frac{1}{p}}(\xi)} \\
&+\theta (\mu,\xi)^{\frac{1}{q}}
\int_{\R^n}\frac{|\eta_2^\mu|^\frac{q+1}{q}}{K^{\frac{1}{q}}(\xi)} \\
& - \theta (\mu,\xi)\irn \langle T \eta^\mu, \eta^\mu \rangle
\bigg],
\end{align*}
and so, since $\theta (\mu,\xi) \eta^\mu \in \Ne_\xi$, in turn we get
\begin{align}
\label{eq:deriv}
\n_\xi I_\xi (\theta (\mu,\xi) \eta^\mu)
= &- \frac{\theta (\mu,\xi)^{\frac{p+1}{p}}}{p+1} \n_\xi Q(\xi)
\int_{\R^n}\frac{|\eta_1^\mu|^\frac{p+1}{p}}{Q^{\frac{p+1}{p}}(\xi)} \\
&- \frac{\theta (\mu,\xi)^{\frac{q+1}{q}}}{q+1} \n_\xi K(\xi)
\int_{\R^n}\frac{|\eta_2^\mu|^\frac{q+1}{q}}{K^{\frac{q+1}{q}}(\xi)}. \notag
\end{align}
From this representation formula, the Mean-Value Theorem and
the local boundedness of $\theta$, the assertion readily follows.
\end{proof}

\subsection{Left and right derivatives of $\Sigma$}
\label{LeR}
Let us define ${\mathbb S}(z)$ as the set of all the positive (dual) solutions $(\eta_1,\eta_2)$ of~\eqref{limit-z}
at the energy level $\Sigma(z)$. The representation formulas for the (left and right) directional
derivatives of $\Sigma$ are provided in the following

\begin{theorem}\label{th:deriv}
The directional derivatives from the left and the
right of $\Sigma$ at every point $z\in\R^n$ along any $w\in\R^n$ exist and it holds
\begin{align*}
&\left(\frac{\partial\Sigma}{\partial w}\right)^{-}\!\!(z)=
\sup_{\eta \in {{\mathbb S}}(z)}\nabla_z I_z (\eta)\cdot w, \\
&\left(\frac{\partial\Sigma}{\partial w}\right)^{+}\!\!(z)=
\inf_{\eta \in {{\mathbb S}}(z)}\nabla_z I_z (\eta)\cdot w.
\end{align*}
Explicitly, we have
\begin{align*}
\left(\frac{\partial\Sigma}{\partial w}\right)^{-}\!\!(z)&=
\sup_{\eta \in {{\mathbb S}}(z)}\left[
-\frac{1}{p+1}\frac{\de Q}{\de w}(z) \int_{\R^n}\frac{|\eta_1|^\frac{p+1}{p}}{Q^{\frac{p+1}{p}}(z)}
-\frac{1}{q+1}\frac{\de K}{\de w}(z) \int_{\R^n}\frac{|\eta_2|^\frac{q+1}{q}}{K^{\frac{q+1}{q}}(z)}
\right]
\\
\left(\frac{\partial\Sigma}{\partial w}\right)^{+}\!\!(z)&=
\inf_{\eta\in {{\mathbb S}}(z)}\left[
-\frac{1}{p+1}\frac{\de Q}{\de w}(z)\int_{\R^n}\frac{|\eta_1|^\frac{p+1}{p}}{Q^{\frac{p+1}{p}}(z)}
-\frac{1}{q+1}\frac{\de K}{\de w}(z)\int_{\R^n}\frac{|\eta_2|^\frac{q+1}{q}}{K^{\frac{q+1}{q}}(z)}
 \right].
\end{align*}
for every $z,w\in\R^n$.
\end{theorem}
\begin{proof}
Let $\{\mu_j\}\subset\R^n$ be a sequence converging to $\mu_0$
and let $\eta^j=\eta^{\mu_j}$ be a sequence of (dual) solutions
of least energy $\Sigma(\mu_j)$. We want to prove that, up to a subsequence,
\begin{equation}
\label{convH-duale}
\eta^{\mu_j}\to \eta^0,\qquad \text{strongly in $\H$},
\qquad   \eta^0\in {\mathbb S}(\mu_0).
\end{equation}
Consider the corresponding (direct) solutions $(u_{\mu_j},v_{\mu_j})$ (resp.\ $(u_0,v_0)$)
of \eqref{eq:Lmu} with $\mu=\mu_j$ (resp.\ $\mu=\mu_0$).
Since $(u_{\mu_j},v_{\mu_j})$ is bounded in $W^{2,(q+1)/q}(\R^n)\times W^{2,(p+1)/p}(\R^n)$
(cf.\ \cite{ASY}), up to a subsequence, it converges weakly to a pair $(u_0,v_0)$. In addition, since
$K$ and $Q$ are uniformly bounded, by virtue of the Schauder local regularity estimates
(cf. \cite{seracta}), $(u_{\mu_j},v_{\mu_j})$ is bounded in $C^{2,\beta}_{\rm loc}(\R^n)$
for some $\beta>0$ and
\begin{equation}
\label{locals}
u_{\mu_j}\to u_0
\quad
\text{and}
\quad
v_{\mu_j}\to v_0,\qquad
\text{locally in $C^2$-sense},
\end{equation}
so that $(u_0,v_0)$ solves \eqref{eq:Lmu} with $\mu=\mu_0$.
We claim that $u_0>0$ and $v_0>0$. By~\cite[Theorem 2]{busi}, for every $j\geq 1$, $u_{\mu_j}$
and $v_{\mu_j}$ are radially symmetric and radially decreasing with
respect to some point, say the origin, that is
\begin{equation}
\label{monotonia}
u_{\mu_j}(x)=u_j(r),\quad v_{\mu_j}(x)=v_j(r),\quad
\frac{d}{dr}u_j(r)<0,\quad\frac{d}{dr}v_j(r)<0,
\end{equation}
for every $r>0$. Hence, for every $j\geq 1$, we have
\begin{align*}
& u_{\mu_j}(0)
\leq -\Delta u_{\mu_j}(0) +u_{\mu_j}(0)
= K(\mu_j) v_{\mu_j}^q(0)\leq \b v_{\mu_j}^q(0),
\\
& v_{\mu_j}(0)
\leq -\Delta v_{\mu_j}(0)+v_{\mu_j}(0)
= Q(\mu_j) u_{\mu_j}^p(0)\leq \b u_{\mu_j}^p(0).
\end{align*}
It follows that, for every $j\geq 1$,
\begin{equation*}
u_{\mu_j}(0)\leq \b^{q+1} u_{\mu_j}^{pq}(0).
\end{equation*}
Then there exists $\hat\delta >0$ such that $u_{\mu_j}(0)\geq\hat\delta$ for every $j\geq 1$.
Similarly, $v_{\mu_j}(0)\geq \hat\delta$ for every $j\geq 1$. Hence, letting
$j\to\infty$, by~\eqref{locals}, we conclude that $u_0(0)\geq \hat\delta$
and $v_0(0)\geq\hat\delta$, which entails $u_0\not\equiv 0$ and $v_0\not\equiv 0$. Since we have
$u_0\geq 0$, $v_0\geq 0$, $K(\mu_0),Q(\mu_0)>0$ and
$$
-\Delta u_0+u_0\geq 0\,\,
\quad
\text{and}
\quad
-\Delta v_0+v_0\geq 0,
$$
the claim just follows by a straightforward application of
the maximum principle.

Observe that, by the continuity of $\S$ and by Fatou's Lemma, we get
$$
\S(\mu_0)=\lim_{j\to\infty}\S(\mu_j)=
\lim_{j\to\infty}I_{\mu_j}(\eta^j)\geq I_{\mu_0}(\eta^0)\geq \S(\mu_0).
$$
Hence
$$
\lim_{j\to\infty}I_{\mu_j}(\eta^j)=I_{\mu_0}(\eta^0)=\S(\mu_0),
$$
which reads as
$$
\lim_{j\to\infty}
\int_{\R^n}\frac{|\eta_1^j|^\frac{p+1}{p}}{Q^{\frac{1}{p}}(\mu_j)}=
\int_{\R^n}\frac{|\eta_1^0|^\frac{p+1}{p}}{Q^{\frac{1}{p}}(\mu_0)},
\qquad
\lim_{j\to\infty}
\int_{\R^n}\frac{|\eta_2^j|^\frac{q+1}{q}}{K^{\frac{1}{q}}(\mu_j)}=
\int_{\R^n}\frac{|\eta_2^0|^\frac{q+1}{q}}{K^{\frac{1}{q}}(\mu_0)}.
$$
In particular, taking into account \eqref{asspoteorig},
for any $\delta>0$, there exists $\rho>0$ such that
$$
\int_{\{|x|\geq\rho\}}|\eta_1^j|^\frac{p+1}{p}<\delta,
\qquad
\int_{\{|x|\geq\rho\}}|\eta_2^j|^\frac{q+1}{q}<\delta,
$$
for every $j\geq 1$ sufficiently large. Moreover, of course
$$
\lim_{j\to\infty}
\int_{\{|x|\leq\rho\}}|\eta_1^j|^\frac{p+1}{p}=
\int_{\{|x|\leq\rho\}}|\eta_1^0|^\frac{p+1}{p},
\qquad
\lim_{j\to\infty}
\int_{\{|x|\leq\rho\}}|\eta_2^j|^\frac{q+1}{q}=
\int_{\{|x|\leq\rho\}}|\eta_2^0|^\frac{q+1}{q}.
$$
Then we have $\eta^{\mu_j}\to \eta^0$ strongly in $\H$,
namely \eqref{convH-duale} holds true.

Without loss of generality, we can prove the formula of the right derivative of $\S$
in the case $n=1$, $z=0$ and $w=1$. For any $\eta^0\in {\mathbb S}(0)$, we get
\begin{align*}
\Sigma(\rho)-\Sigma(0)&\leq I_\rho(\vartheta(\rho,0)\eta^0)-I_0(\eta^0)
\\
&=\rho \n_\xi I_\xi (\vartheta (\xi,0) \eta^0)|_{\xi=\mu \in[0,\rho]}.
\end{align*}
Whence, by virtue of \eqref{eq:deriv} and the arbitrariness of $\eta^0 \in {\mathbb S}(0)$,
\begin{equation*}
\limsup_{\rho\to 0^+}\frac{\Sigma(\rho)-\Sigma(0)}{\rho}
\leq \inf_{\eta^0\in {\mathbb S}(0)}\bigg[
- \frac{Q'(0)}{p+1} \int_{\R^n}\frac{|\eta_1^0|^\frac{p+1}{p}}{Q^{\frac{p+1}{p}}(0)}
- \frac{K'(0)}{q+1}  \int_{\R^n}\frac{|\eta_2^0|^\frac{q+1}{q}}{K^{\frac{q+1}{q}}(0)} \bigg].
\end{equation*}
Moreover, similarly, we get
\begin{align*}
\Sigma(\rho)-\Sigma(0)&\geq
I_\rho(\vartheta(\rho,\rho) \eta^\rho)
-I_0(\vartheta(0,\rho)\eta^\rho) \\
&=\rho\n_\xi I_\xi (\vartheta (\xi,\rho)\eta^\rho)|_{\xi=\mu \in[0,\rho]},
\end{align*}
so that, by exploiting \eqref{eq:deriv} and \eqref{convH-duale}, we conclude
\begin{equation*}
\liminf_{\rho\to 0^+}\frac{\Sigma(\rho)-\Sigma(0)}{\rho}
\geq \inf_{\eta^0\in {\mathbb S}(0)}\bigg[
- \frac{Q'(0)}{p+1} \int_{\R^n}\frac{|\eta_1^0|^\frac{p+1}{p}}{Q^{\frac{p+1}{p}}(0)}
- \frac{K'(0)}{q+1}  \int_{\R^n}\frac{|\eta_2^0|^\frac{q+1}{q}}{K^{\frac{q+1}{q}}(0)} \bigg].
\end{equation*}
Then the desired formula for the right derivative of $\Sigma$ follows.
A very similar argument provides the corresponding formula for the left derivative.
\end{proof}

\begin{remark}\rm
\label{smoothre}
Nowadays, the further regularity of $\Sigma$ is, to our knowledge, an open problem. Actually, not
even in the case of a single equation the situation is very well
understood. For instance, on one hand, if we consider the problem
\begin{equation*}
-\eps^2 \Delta u+V(x)u=K(x)u^p \quad \text{in $\R^n$},
\qquad u>0\quad \text{in $\R^n$},
\end{equation*}
then $\Sigma\in C^m(\R^n)$ provided that both the potentials
$V$ and $K$ belong to $C^m(\R^n)$, with $m\geq 1$. On the other hand,
if $f$ is not a power (and does not satisfy conditions
ensuring uniqueness up to translations), for the equation
\begin{equation*}
-\eps^2\Delta u+V(x)u=K(x)f(u) \quad \text{in $\R^n$},
\qquad u>0\quad \text{in $\R^n$},
\end{equation*}
we do not know which regularity beyond ${\rm Lip}_{\rm loc}$ can be achieved by $\S$.
Even though we do not possess any specific
counterexample, our feeling is that there exist
functions $f$ for which the associated $\Sigma$ fails to be $C^1$ smooth. It is evident by the
(left and right) derivative formulas of $\Sigma$ that
its further regularity is related to the uniqueness of
positive radial solutions to $-\Delta u+u=f(u)$ in $\R^n$,
$u>0$ in $\R^n$, which occurs just for very particular nonlinearities $f$. Based
upon these considerations, for semilinear systems, the further
regularity of $\Sigma$ seems an {\em ever harder} matter, since as already stressed nothing
is known, so far, about the uniqueness of solutions to the system
\begin{equation*}
-\Delta u+u=f(v), \qquad
-\Delta v+v=g(u), \quad\text{in $\R^n$},\qquad u,v>0\quad \text{in $\R^n$},
\end{equation*}
not even with the particular choices $f(v)=v^q$ and $g(u)=u^p$.
\end{remark}

\section{\bf Proof of the results}
\label{PROOFS}
\def\theequation{3.\arabic{equation}}\makeatother
\setcounter{equation}{0}

\subsection{Proof of Theorem~\ref{neceMAIN}}
Let $z\in{\mathcal E}$ and let $(u_{\eps_h},v_{\eps_h})\in
W^{2,\frac{q+1}{q}}(\R^n)\times W^{2,\frac{p+1}{p}}(\R^n)$ be a corresponding
a sequence of strong solutions to~\eqref{problema} with $|u_{\eps_h}(z)|,|v_{\eps_h}(z)|\geq\delta$
for some $\delta>0$, $|u_{\eps_h}(z+\eps_hx)|\to 0,|v_{\eps_h}(z+\eps_hx)|\to 0$
as $|x|\to\infty$ uniformly w.r.t.\ $h$, and
${\eps_h}^{-n}f_{\eps_h}(u_{\eps_h},v_{\eps_h})
\to\Sigma(z)$ as $h\to\infty$. Let us set:
$$
\vfi_h (x)=u_{\eps_h}(z+\eps_h x)
\qquad
\text{and}
\qquad
\psi_h (x)=v_{\eps_h}(z+\eps_h x),
$$
for all $h\geq 1$. Then, since $(u_{\eps_h},v_{\eps_h})$ is a solution
\eqref{problema}, $(\vfi_h,\psi_h)$ is solution of
\begin{equation}
\label{Sh}
-\Delta \vfi_h+\vfi_h=K(z+\eps_h x)\psi_h^q,  \qquad
-\Delta \psi_h+\psi_h=Q(z+\eps_h x)\vfi_h^p.
\end{equation}
By arguing as in the proof of Theorem \ref{th:deriv},
it is readily proved that, up to a subsequence, $(\vfi_h)$ and $(\psi_h)$
converge weakly in $W^{2,(q+1)/q}(\R^n)\times W^{2,(p+1)/p}(\R^n)$
to some $\varphi_0$ and $\psi_0$ respectively.
Let us now prove that there exist $\Theta>0$,
$\rho>0$ and $h_0\geq 1$ such that
\begin{equation}
\label{decayll}
\vfi_{h}(x)\leq ce^{-\Theta|x|}
\quad
\text{and}
\quad
\psi_{h}(x)\leq ce^{-\Theta|x|},
\qquad\text{for all $|x|\geq \rho$ and $h\geq h_0$}.
\end{equation}
We follow the line of \cite{dFJ}.
Since $z\in{\mathcal E}$, then the functions $\vfi_{h}$ and $\psi_{h}$ decay
to zero at infinity, uniformly with respect to $h$. Hence, since $p,q>1$,
we can find $\rho>0$, $\Theta>0$ and $h_0\geq 1$ such that
\begin{align*}
K(z+\eps_h x)\psi_h^{q} & \leq (1-\Theta^2)\psi_h, \\
\noalign{\vskip3pt}
Q(z+\eps_h x)\vfi_h^{p} & \leq (1-\Theta^2)\vfi_h,
\end{align*}
for all $|x|>\rho$ and $h\geq h_0$. Let us set
$$
\xi(x)=\mu e^{-\Theta(|x|-\rho)},\qquad
\mu=\max_{|x|=\rho}\max_{h\geq h_0} (\psi_h+\vfi_h),
$$
and introduce the set
$$
A=\bigcup_{R>\rho}D_R,
$$
where, for any $R>\rho$, we put
$$
\quad D_R=\big\{\rho<|x|<R:\,\,
\psi_h(x)+\vfi_h(x)>\xi(x)\,\,\,\,
\text{for some $h\geq h_0$}\big\}.
$$
If $A=\emptyset$, we are done. Instead, if $A$
is nonempty, there exists $R_*>\rho$ such that
\begin{align*}
\Delta(\xi-\psi_h-\vfi_h) & \leq\Big[\Theta^2-\frac{\Theta(n-1)}{|x|}\Big]\xi(x)
-\Theta^2 \psi_h-\Theta^2 \vfi_h \\
& \leq\Theta^2 (\xi-\psi_h-\vfi_h)<0,
\qquad\text{on $D_R$ for all $R\geq R_*$}.
\end{align*}
Hence, by the maximum principle, since $(\xi-\psi_h-\vfi_h)|_{\{|x|=\rho\}}\geq 0$, we get
$$
\xi-\psi_h-\vfi_h\geq \min\big\{0,\min_{|x|=R}(\xi-\psi_h-\vfi_h)\big\},
\qquad\text{for all $R\geq R_*$}
$$
so that, letting $R\to\infty$, yields, for any $\rho>0$,
$\psi_h(x)+\vfi_h(x)\leq \xi(x)$ for $|x|>\rho$,
which contradicts the definition of $D_{R_*}\not=\emptyset$.

By virtue of the Schauder interior estimates (see e.g.\ \cite{seracta}), $\vfi_{h}\to \vfi_0$ and
$\psi_{h}\to \psi_0$ locally in $C^2$ sense, so that $(\vfi_0,\psi_0)$
is a (nontrivial, radial, decaying) solution to \eqref{limit-z}.
Moreover, in light of the exponential barriers provided by \eqref{decayll},
since $z \in \mathcal{E}$, it is not difficult to see that
$(\vfi_0,\psi_0) \in {\mathbb S}(z)$, for we have
\begin{align*}
\S(z)&= {\textstyle\big(\frac 12 - \frac{1}{q+1}\big)} \irn K(z)|\psi_0|^{q+1}
+{\textstyle\big(\frac 12 - \frac{1}{p+1}\big)} \irn Q(z)|\vfi_0|^{p+1} \\
&=f_z(\vfi_0,\psi_0)=I_z(\eta^0),
\end{align*}
where $\eta^0$ is the dual solution corresponding to $(\vfi_0,\psi_0)$.

Let us now consider the Lagrangian $\mathscr{L}: \Rn \times \R \times
\R \times \Rn \times \Rn \to
\R$ defined as
\[
\mathscr{L}(x,s_1,s_2,\xi_1,\xi_2)=
\xi_1 \cdot \xi_2 +s_1 s_2
-{\textstyle\frac{1}{q+1}}K(z+\eps_h x)s_2^{q+1}
-{\textstyle\frac{1}{p+1}}Q(z+\eps_h x)s_1^{p+1}.
\]
Then system \eqref{Sh} rewrites as
\[
\begin{cases}
-{\rm div}\left(\partial_{\xi_2}\mathscr{L}(x,\vfi_h,\psi_h,\n \vfi_h,\n \psi_h) \right)
+ \partial_{s_2}\mathscr{L}(x,\vfi_h,\psi_h,\n \vfi_h,\n \psi_h)=0, &  \,\,\text{in $\R^n$},
\\
-{\rm div}\left(\partial_{\xi_1}\mathscr{L}(x,\vfi_h,\psi_h,\n \vfi_h,\n \psi_h) \right)
+ \partial_{s_1}\mathscr{L}(x,\vfi_h,\psi_h,\n \vfi_h,\n \psi_h)=0, &  \,\,\text{in $\R^n$}, \\
\noalign{\vskip2pt}
\, \vfi_h,\,\psi_h>0, &  \,\,\text{in $\R^n$}.
\end{cases}
\]
By the Pucci-Serrin identity for systems \cite[see \textsection 5]{puc-ser}, we have
\begin{gather*}
\sum\limits^n_{i,\,l=1}
\irn \partial_i {\boldsymbol q}^l
\partial_{(\xi_2)_i}\mathscr{L}(x,\vfi_h,\psi_h,\n \vfi_h,\n \psi_h) \partial_l \psi_h
\\
+ \sum\limits^n_{i,\,l=1}
\irn \partial_i {\boldsymbol q}^l
\partial_{(\xi_1)_i}\mathscr{L}(x,\vfi_h,\psi_h,\n \vfi_h,\n \psi_h) \partial_l \vfi_h
\\
= \irn \big[({\rm div}\, {\boldsymbol q})\mathscr{L}(x,\vfi_h,\psi_h,\n \vfi_h,\n \psi_h)
+ {\boldsymbol q} \cdot \partial_x  \mathscr{L}(x,\vfi_h,\psi_h,\n \vfi_h,\n \psi_h)     \big],
\end{gather*}
for all ${\boldsymbol q} \in C^1_{\rm c}\left(\Rn, \Rn \right)$.
Let us take, for $\lambda>0$,
$$
{\boldsymbol q}(x) =  (\Upsilon(\lambda x), \,0, \ldots, \,0),
$$
and $\Upsilon\in C^1_{{\rm c}} (\Rn)$
such that $\Upsilon(x)=1$ if $|x|\leq 1$ and $\Upsilon(x)=0$
if $|x|\geq 2$. Then,
\begin{gather*}
\sum\limits^n_{i=1}
\irn \lambda \partial_i \Upsilon(\lambda x) \partial_i \vfi_h \, \partial_1 \psi_h
+\sum\limits^n_{i=1}
\irn \lambda \partial_i \Upsilon(\lambda x) \partial_i \psi_h \, \partial_1 \vfi_h
\\
=\irn \lambda \partial_1 \Upsilon(\lambda x) \mathscr{L}(x,\vfi_h,\psi_h,\n \vfi_h,\n \psi_h)
\\
+\irn \eps_h \Upsilon(\lambda x)
\left[-{\textstyle\frac{1}{q+1}}\partial_1 K(z+\eps_h x)\psi_h^{q+1}
-{\textstyle\frac{1}{p+1}}\partial_1 Q(z+\eps_h x)\vfi_h^{p+1} \right].
\end{gather*}
By the arbitrariness of $\lambda>0$, letting
$\lambda \to 0$ and keeping $h$ fixed, we obtain
\[
\irn \left[
-{\textstyle\frac{1}{q+1}}\partial_1 K(z+\eps_h x)\psi_h^{q+1}
-{\textstyle\frac{1}{p+1}}\partial_1 Q(z+\eps_h x)\vfi_h^{p+1} \right]=0.
\]
Therefore, letting now $h\to \infty$, since in light of~\eqref{assp2} we get
\begin{equation*}
|\nabla K(z+\eps_hx)|,
|\nabla Q(z+\eps_hx)|\leq ce^{M\eps_h|x|},
\qquad\text{for $|x|$ large,}
\end{equation*}
by virtue of \eqref{decayll}, there holds
\[
\irn \left[
-{\textstyle\frac{1}{q+1}} \partial_1 K(z)\psi_0^{q+1}
-{\textstyle\frac{1}{p+1}} \partial_1 Q(z)\vfi_0^{p+1} \right]=0.
\]
Analogously, we can show that, for all $w \in \Rn$,
\[
\irn \left[
-{\textstyle\frac{1}{q+1}} \n K(z)\psi_0^{q+1}
-{\textstyle\frac{1}{p+1}} \n Q(z)\vfi_0^{p+1} \right] \cdot w=0.
\]
Hence
\begin{equation}
\label{form-finale}
-{\textstyle\frac{1}{p+1}}\frac{\de Q}{\de w}(z) \int_{\R^n}\frac{|\eta^0_1|^\frac{p+1}{p}}{Q^{\frac{p+1}{p}}(z)}
-{\textstyle\frac{1}{q+1}}\frac{\de K}{\de w}(z) \int_{\R^n}\frac{|\eta^0_2|^\frac{q+1}{q}}{K^{\frac{q+1}{q}}(z)}
=0
\end{equation}
Since $\eta^0\in{\mathbb S}(z)$, by Theorem \ref{th:deriv} we have
\begin{align*}
\left(\frac{\partial\Sigma}{\partial w}\right)^{+}\!\!(z)&=
\inf_{\eta\in {\mathbb S}(z)}\nabla_z I_z(\eta)\cdot w  \\
&\leq
-{\textstyle\frac{1}{p+1}}\frac{\de Q}{\de w}(z) \int_{\R^n}\frac{|\eta^0_1|^\frac{p+1}{p}}{Q^{\frac{p+1}{p}}(z)}
-{\textstyle\frac{1}{q+1}}\frac{\de K}{\de w}(z) \int_{\R^n}\frac{|\eta^0_2|^\frac{q+1}{q}}{K^{\frac{q+1}{q}}(z)}
=0.
\end{align*}
Then, by the very definition of $(-\Sigma)^0(z;w)$ (see Definition~\ref{def:der}), we get
\begin{equation*}
(-\Sigma)^0(z;w)\geq\left(\frac{\partial (-\Sigma)}
{\partial w}\right)^{+}\!\!(z)\geq 0,\quad
\text{for every $w\in\R^n$}.
\end{equation*}
Then $0\in \partial_C(-\Sigma)(z)$ and, since
$\partial_C (-\Sigma)(z)=-\partial_C \Sigma(z)$ (cf.\ \cite{clarke}), we
obtain $z\in\mathcal{K}$.

\subsection{Proof of Corollary~\ref{uniTH}} It suffices to combine
Theorems~\ref{neceMAIN} and~\ref{th:deriv}, taking into account
what discussed in Section~\ref{uniqrem} about the conjectured
explicit representation formula for $\S$.

\subsection{Proof of Theorem~\ref{neceMAIN2}} Let $m\geq 1$ and $z\in\mathcal{E}_m$. The assertion
follows by mimicking the various steps in the proof of
Theorem~\ref{neceMAIN} with $\mathcal{E}_m$ in place of
$\mathcal{E}$, and combining formula \eqref{form-finale}
with the definitions of $\Gamma_{z,m}^\mp$ and $\mathcal{K}_m$,
taking into account that $\eta^0\in\GG_m(z)$, as it holds $I_z(\eta^0)=m$,
being $\eta^0$ the strong limit of $\eta^{\eps_j}$.
Indeed, by \eqref{form-finale}, there holds
$$
\Gamma_{z,m}^+(w)\geq 0,\quad\forall w\in\R^n,
\qquad
\Gamma_{z,m}^-(w)\geq 0,\quad\forall w\in\R^n,
$$
so that $0\in\partial\Gamma_{z,m}^+(0)\cap\partial\Gamma_{z,m}^-(0)$,
yielding $z\in\mathcal{K}_m$.

\end{document}